\documentclass[]{interact}

\usepackage{epstopdf}
\usepackage[caption=false]{subfig}

\usepackage[numbers,sort&compress]{natbib}

\makeatletter
\def\NAT@def@citea{\def\@citea{\NAT@separator}}
\makeatother
\usepackage{anyfontsize}
\theoremstyle{plain}
\newtheorem{theorem}{Theorem}[section]
\newtheorem{lemma}[theorem]{Lemma}

\newtheorem{proposition}[theorem]{Proposition}

\theoremstyle{definition}
\newtheorem{definition}[theorem]{Definition}

\usepackage{xcolor}

\theoremstyle{remark}

\begin{document}


\title{Riesz Bases in Krein Spaces}
\author{
	\name{Shah Jahan\textsuperscript{a},  P. Sam  Johnson\textsuperscript{b}}
	\affil{\textsuperscript{a} Department of Mathematics, Central University  of Haryana, Mahendergarh 123029, Haryana, India. Email address: shahjahan@cuh.ac.in} \textsuperscript{b} Department of Mathematical and Computational Sciences,
	National Institute of Technology Karnataka (NITK), Surathkal, Mangaluru 575025, India. Email address: sam@nitk.edu.in}

\maketitle

\begin{abstract}
We start by introducing and studying the definition of a Riesz basis in a Krein space $(\mathcal{K},[.,.])$, along with a condition under which a Riesz basis becomes a Bessel sequence. The concept of biorthogonal sequence in Krein spaces is also introduced, providing an equivalent characterization of a Riesz basis. Additionally, we explore the concept of the Gram matrix, defined as the sum of a positive and a negative Gram matrices, and specify conditions under which the Gram matrix becomes bounded in Krein spaces. Further, we characterize the conditions under which the Gram matrices $\{[f_n,f_j]_{n,j \in I_+}\}$ and $\{[f_n,f_j]_{n,j \in I_-}\}$  become bounded invertible operators. Finally, we provide an equivalent characterization of a Riesz basis in terms of Gram matrices.
\end{abstract}

\begin{amscode}42C15, 46C05, 46C20.\end{amscode}

\maketitle

\section{Introduction}\label{sec1}
Hilbert space frames were originally introduced by Duffin and Schaeffer \cite{g} in 1952 to deal with some problems in non-harmonic Fourier analysis. After some decades, Daubechies, Grossmann and Meyer \cite{f} announced formally the definition of frame in the abstract Hilbert spaces in $1986$.  After their work, frame theory began to be widely used, particularly in the more specialized context of wavelet frames and Gabor frames. The linear independence property for a basis, which allows every vector to be uniquely represented as a linear combination, is very restrictive for practical problems. Frames allow each element in the space to be written as a linear combination of the elements in the frame, but linear independence is not required. Frames can be viewed as redundant bases which are generalizations of Riesz bases. This redundancy property  is extremely important in applications.

	Let $H$ be a Hilbert space and $I$ be a countable index set. A collection $\{f_n\}_{n \in I}$ in $H$ is called a frame for $H$ if there exist positive constants
	$0< A\leqslant B< \infty$ such that
	\begin{align}\label{h1}
	A\|f\|^2 \leqslant \sum_{n \in I} |\langle f, f_n \rangle|^2 \leqslant B\|f\|^2 ,\quad \text{for all} \ f \in H. \end{align} The bounded linear operator $S:H\rightarrow H$ defined by \begin{align*}
	Sf=\sum_{n \in I}\langle f,f_n \rangle f_n,~~~~ f \in H
	\end{align*} is known as the frame operator associated to the frame $\{f_n\}_{n \in I}.$ This operator $S$ is bounded invertible, positive and self adjoint. It allows to reconstruct each vector in terms of the sequence $\{f_n\}_{n \in I}$ as follows: \begin{align}\label{h2} f = \sum_{n \in I}\langle f,S^{-1}f_n \rangle f_n=\sum_{n \in I}\langle f,f_n \rangle S^{-1}f_n.\end{align} The formula (\ref{h2}) is known as reconstruction formula associated to $\{f_n\}_{n \in I}$. If $S=I,$ then the reconstruction formula resembles the Fourier series of $f$ associated with the orthonormal sequence $\{f_n\}_{n \in I}.$\\ 
	
	Giribet et al. \cite{k} introduced and studied frames for Krein spaces. The authors have given a new notion of frames in Krein spaces in \cite{s}. In this paper, we study the concept of Riesz basis in Krein spaces by decomposing basis in a natural way and obtain some new results on frames and Riesz bases. We organize the paper as follows. In Section 2, definition and basic results  of a Krein space are given which are required in this paper. In Section 3, we define the concept of a Riesz basis in Krein spaces and a necessary and sufficient condition is given for a sequence to be a Riesz basis. In Section 4, we define Gram matrix in Krein spaces and study  properties of Gram matrix.  Further, we show that if $\{f_n\}_{n \in I}$ is a Bessel sequence, then the Gram matrix defines bounded invertible operators on $\ell_2(I_+)$ and $\ell_2(I_-).$ 
	\section{Preliminaries} We begin this section with the definition of Krein space and some notations which are used in sequel. For more detail about Krein spaces, one may refer to \cite{a,b,d}. 
	\begin{definition}A linear space $\mathcal{K}$ with a $Q$-metric $[x,y],$ which admits a canonical decomposition $\mathcal{K}=\mathcal{K}^{+} \oplus \mathcal{K}^{-}$ in which $\mathcal{K}^{+}$ and $\mathcal{K}^{-}$ are complete relative to the norm $\|x\|=[x,x]^\frac{1}{2}$ and $\|x\|=(-[x,x])^\frac{1}{2}$ respectively is called a Krein space.
	\end{definition}Through out this paper $(\mathcal{K},[.,.])$ denotes a Krein space, which means that $\mathcal{K}$ can be written as the direct $[.,.]$-orthogonal sum $\mathcal{K}^+\oplus \mathcal{K}^- $ of Hilbert spaces.
	Two vectors in  $ \mathcal{K}$ are said to be orthogonal if $[x,y]=0$ and we denote this by the symbol $x[\bot]y.$ By orthonormal basis, we mean a family $\{e_i\}_{i \in I} \in \mathcal{K}$ such that $[e_i,e_j]=\pm\delta_{ij},$ for all $i,j \in I.$ If $I_+=\{n \in I : [f_n,f_n]\geq 0\}$ and $I_-=\{n \in I : [f_n,f_n]< 0\},$ then the orthogonal decomposition of $\ell_2(I)$ is given by $\ell_2(I)=\ell_2(I_+)\oplus \ell_2(I_-).$
	
	\begin{lemma}\label{rl}
		Let $( \mathcal{K},[.,.])$ be a Krein space.
		\begin{enumerate}
			\item [(a)] If $x,y \in \mathcal{K^+}$ satisfying $[x,z]=[y,z],  \text{for all}  z \in \mathcal{K^+}$, then $x=y$.
			\item [(b)] If $x,y \in \mathcal{K^-}$ satisfying $[x,z]=[y,z], \text{for all}  z \in \mathcal{K^-}$, then $x=y$.
		\end{enumerate}
	\end{lemma}

	\section{Riesz basis} We begin this section with the definition of Riesz basis in Krein spaces  and followed by we prove that a Riesz basis is a Bessel sequence. 
	\begin{definition} A Riesz basis for  $( \mathcal{K},[.,.])$ is a system of the form $\mathcal{F}= (\{U^+e_n\}_{n \in I_+}\cup\{U^-e_n\}_{n \in I_-})$, where $\{e_n\}_{n \in I_+}$ is an orthonormal basis for $\mathcal{K}^+$ and $\{e_n\}_{n \in I_-}$ is an orthonormal basis for $\mathcal{K}^-$ and $U^+: \mathcal{K^+} \rightarrow \mathcal{K^+}$, $U^-: \mathcal{K^-} \rightarrow \mathcal{K^-}$ are bounded bijective operators.  
	\end{definition}
	\begin{theorem}\label{T1}
		Let $\{f_n\}_{n \in I}$ be a Riesz basis for  $( \mathcal{K},[.,.]).$ Then $\{f_n\}_{n \in I}$ is a Bessel sequence with Bessel bounds $\|U^+\|$ and $\|U^-\|$ respectively. Moreover, there exist unique sequences $\{g_n\}_{n \in I_+} \in \mathcal{K^+}$ and $\{g_n\}_{n \in I_-} \in \mathcal{K^-}$ such that
		\begin{equation}
		f = \sum_{n \in I_+}[f,g_n]f_n, \ \mbox{for all} \ f \in \mathcal{K^+}
		\end{equation}and
		\begin{equation}
		f = \sum_{n \in I_-}[f,g_n]f_n, \ \mbox{for all} \ f \in \mathcal{K^-}.
		\end{equation} The sequences $\{g_n\}_{n \in I_+}$ and $\{g_n\}_{n \in I_-}$ are Riesz bases and both the series converge unconditionally.
	\end{theorem}
	\begin{proof} Since $\{f_n\}_{n \in I}$ is a Riesz basis for the Krein space $( \mathcal{K},[.,.])$, we can write $\{f_n\}_{n \in I}= (\{U^+e_n\}_{n \in I_+}\cup\{U^-e_n\}_{n \in I_-})$, where $\{e_n\}_{n \in I_+}$ is an orthonormal basis for $\mathcal{K}^+$,  $\{e_n\}_{n \in I_-}$ is an orthonormal basis for $\mathcal{K}^-$ and  $U^+: \mathcal{K^+} \rightarrow \mathcal{K^+}$, $U^-: \mathcal{K^-} \rightarrow \mathcal{K^-}$ are bounded bijective operators. Let $f \in \mathcal{K^+}$. By expanding $(U^+)^{-1}$ in term of orthonormal basis $\{e_n\}_{n \in I_+}$, we have 
		\begin{eqnarray*}
			(U^+)^{-1}f&=&\sum_{n \in I_+}[(U^+)^{-1}f,e_n]e_n
			\\&=&\sum_{n \in I_+}[f, ((U^+)^{-1})^*e_n]e_n.
		\end{eqnarray*}
		Now with $g_n=((U^+)^{-1})^*e_n,$ we can write
		\begin{eqnarray*}
			f&=&U^+(U^+)^{-1}f 
			\\&=&\sum_{n \in I_+}[f, ((U^+)^{-1})^*e_n]U^+e_n
			\\&=&\sum_{n \in I_+}[f,g_n]f_n.
		\end{eqnarray*}And similarly, one can easily prove that for all $f \in \mathcal{K^-}$ and with $g_n=((U^-)^{-1})^*e_n,$ we have 
		\begin{eqnarray*}
			f &=&\sum_{n \in I_-}[f, ((U^-)^{-1})^*e_n]U^-e_n
			\\&=&\sum_{n \in I_-}[f,g_n]f_n.
		\end{eqnarray*}
		Since both the operators $((U^+)^{-1})^*$ and $((U^-)^{-1})^*$ are bounded and bijective, so $\{g_n\}_{n \in I_+}$ and $\{g_n\}_{n \in I_-}$ are Riesz bases for $\mathcal{K^+}$ and $\mathcal{K^-}$ respectively. Now for $f \in \mathcal{K^+}$, we have
		\begin{eqnarray*}\sum_{n \in I_+}|[f,f_n]|^2 &=& \sum_{n \in I_+}|[f,U^+e_n]|^2
			\\&=&\|(U^+)^*f\|^2
			\\&\leq& \|(U^+)^*\|\|f\|
			\\&=& \|U^+\|\|f\|.\end{eqnarray*}
		So $\{f_n\}_{n \in I_+}$ is a Bessel sequence with Bessel bound $\|U^+\|$ and for $f \in \mathcal{K^-}$ one can easily prove that $\{f_n\}_{n \in I_-}$ is also a Bessel sequence with Bessel bound $\|U^-\|.$ Hence $\{f_n\}_{n \in I}$ is a Bessel sequence for the Krein space $( \mathcal{K},[.,.]).$ The proof for the unconditional convergence of the series follows from the Corollary 3.3 \cite{s}. Now, to complete the proof we only need to show  $\{g_n\}_{n \in I_+}$ and $\{g_n\}_{n \in I_-}$ are unique sequences in  $\mathcal{K^+}$ and $\mathcal{K^-}$ respectively. For this, suppose that $\{g_n\}_{n \in I_+}$ and $\{h_n\}_{n \in I_+}$ are two sequences in $\mathcal{K^+}$ such that 
		\begin{eqnarray}\label{r1}
		f= \sum_{n \in I_+}[f,g_n]h_n = \sum_{n \in I_+}[f,h_n]g_n, \ \mbox{for all} \ f \in \mathcal{K^+}.
		\end{eqnarray}  Then $g_n = h_n$ for all $f \in \mathcal{K^+}.$ Hence from (\ref{r1}), we  have
		\begin{equation*}
		[f,g_n] = [f,h_n], \  \mbox{for all} \ f \in \mathcal{K^+}.
		\end{equation*}
		Similarly, suppose that $\{g_n\}_{n \in I_-}$ and  $\{h_n\}_{n \in I_-}$ are two sequences in $\mathcal{K^-}$, then $[f,g_n] = [f,h_n], \mbox{ for all } f \in \mathcal{K^-}.$ From these equalities and by Lemma \ref{rl}, the result for the uniqueness follows.
	\end{proof}
\noindent	The following is the definition of biorthogonal sequence in a Krein space $(\mathcal{K},[.,.]).$
	\begin{definition} Two sequences $\{f_n\}_{n \in I}$ and $\{g_n\}_{n \in I}$ are said to be biorthogonal in a Krein space $( \mathcal{K},[.,.])$ if 
		\begin{equation}[f_n, g_m]= \pm \delta_{nm}. \end{equation}
	\end{definition}
	
	\begin{theorem} Let $\{f_n\}_{n \in I}= (\{U^+e_n\}_{n \in I_+}\cup\{U^-e_n\}_{n \in I_-})$ be a Riesz basis for  $( \mathcal{K},[.,.])$.  Then there exist positive constants $A,B,A^{'} ,B^{'}$ such that
		\begin{equation}
		A\|f\|^2 \leq \sum_{n \in I_+}|[f,f_n]|^2\leq B\|f\|^2, \ \mbox{for all} \ f \in \mathcal{K^+}
		\end{equation} and 
		\begin{equation}
		A^{'}\|f\|^2 \leq \sum_{n \in I_-}|[f,f_n]|^2\leq B^{'}\|f\|^2, \ \mbox{for all} \ f \in \mathcal{K^-}.
		\end{equation}
		Moreover, the largest and the smallest values for these constants $A,A^{'}, B,B^{'}$ are $\frac{1}{\|(U^+)^{-1}\|^2}$, $\frac{1}{\|(U^-)^{-1}\|^2}$, $\|U^+\|^2$ and $\|U^{-1}\|^2$ respectively.
	\end{theorem}
	\begin{proof}
		By Theorem \ref{T1}, if $\{f_n\}_{n \in I}=(\{U^+e_n\}_{n \in I_+}\cup\{U^-e_n\}_{n \in I_-})$ is a Riesz basis for the Krein space $( \mathcal{K},[.,.]),$ then it is a Bessel sequence, with bounds $\|U^+\|^2$ and $\|U^{-1}\|.$ Next we prove the result for the lower bounds. For this, let $f \in \mathcal{K^+},$
	\begin{equation}
		\sum_{n \in I_+}|[f, f_n]|^2 = \sum_{n \in I_+}|[f, U^+e_n]| = \|(U^+)^*f\|^2
		\label{O1}
	\end{equation}
		and we can write 
		\begin{eqnarray}
		\notag \|f\|&=&\|((U^+)^*)^{-1}(U^+)^*f\|
		\\&\leq&\notag \|((U^+)^*)^{-1}\| \|(U^+)^*f\|
		\label{O2}\\&=& \|(U^+)^{-1}\| \|(U^+)^*f\|.
		\end{eqnarray}
		Therefore from (\ref{O1}) and (\ref{O2}), we conclude that 
		\begin{eqnarray*}
			\frac{1}{\|(U^+)^{-1}\|^2}\|f\|^2 \leq \sum_{n \in I_+}|[f, f_n]|^2, \ \mbox{for all} \ f \in \mathcal{K^+},
		\end{eqnarray*} and similarly for $f \in \mathcal{K^-},$ one can easily prove that 
		\begin{eqnarray*}
			\frac{1}{\|(U^-)^{-1}\|^2}\|f\|^2 \leq \sum_{n \in I_-}|[f, f_n]|^2, \ \mbox{for all} f \in \mathcal{K^-}.
		\end{eqnarray*}
	\end{proof}

Let $\mathcal{K}_1$ and $\mathcal{K}_2$ be  Krein spaces such that $\mathcal{K}_1 = \mathcal{K}^+_1\oplus \mathcal{K}^+_1$ and $\mathcal{K}_2 = \mathcal{K^+}_2\oplus \mathcal{K^-}_2.$ Let $\{h_n\}_{n \in I_+}, \{g_n\}_{n \in I_+}, \{h_n\}_{n \in I_-}$ and $\{g_n\}_{n \in I_-}$ be  sequences in $\mathcal{K}^+_1$, $\mathcal{K}^+_2$,  $\mathcal{K}^-_1$ and $\mathcal{K}^-_2$ respectively. 
If the sequences $\{g_n\}_{n \in I_+}$ and $\{g_n\}_{n \in I_-}$ are Bessel sequences with bounds $B$ and $B^{'}$ respectively and if  $\{h_n\}_{n \in I_+}, \{h_n\}_{n \in I_-}$ are complete in $\mathcal{K}^+_1$ and $\mathcal{K}^-_2$ respectively,  there exist positive constants $A$ and $A^{'}$ such that 
		\begin{eqnarray*}
			A\sum_{n \in I_+}|c_n|^2 \leq \Big\|\sum_{n \in I_+}c_nh_n\Big\|^2, \ \text{for all scalar sequences }\{c_n\}_{n \in I_+} \in \ell_2(I_+)
		\end{eqnarray*} and \begin{eqnarray*}
			A^{'}\sum_{n \in I_-}|c_n|^2 \leq \Big\|\sum_{n \in I_-}c_nh_n\Big\|^2,\  \text{for all scalar sequences }\{c_n\}_{n \in I_-} \in \ell_2(I_-).
		\end{eqnarray*}

	Since $\{h_n\}_{n \in I_+}$ is complete in $\mathcal{K}^+_1$, every $h \in \text{span}\{h_n ; n \in I_+\}$ has a unique representation $h= \sum_{n \in I_+}c_nh_n$, where $\{c_n\}_{n \in I_+} \in \ell_2(I_+).$ 
	Then \begin{eqnarray*}
		U^+\Big(\sum_{n \in I_+}c_nh_n\Big)= \sum_{n \in I_+}c_ng_n, \ \text{for all scalar  sequences }\{c_n\}_{n \in I_+} \in \ell_2(I_+)
	\end{eqnarray*} and 
		\begin{eqnarray*}
			U^-\Big(\sum_{n \in I_-}c_nh_n\Big)= \sum_{n \in I_-}c_ng_n, \ \text{for all scalar sequences }\{c_n\}_{n \in I_-} \in \ell_2(I_-),\end{eqnarray*} define bounded linear operators from $\text{span}\{h_n: n \in I_+\}$ into $\text{span}\{g_n: n \in I_+\}$ and $\text{span}\{h_n: n \in I_-\}$ into $\text{span}\{g_n: n \in I_-\}$ respectively.  
		 So $U^+$ is a well-defined linear operator and for $\{c_n\}_{n \in I_+} \in \ell_2(I_+)$ we have, 
		\begin{eqnarray*}
			\Big\|U^+\Big(\sum_{n \in I_+}c_nh_n\Big)\Big\|^2&=&\Big\|\sum_{n \in I_+}c_ng_n \Big\|^2
			\\&\leq&B\sum_{n \in I_+}|c_n|^2
			\\&\leq&\frac{B}{A}\sum_{n \in I_+}|c_n|^2.
		\end{eqnarray*}Therefore,
		\begin{eqnarray*}
			\Big\|U^+\Big(\sum_{n \in I_+}c_nh_n\Big)\Big\| &\leq&\sqrt{\frac{B}{A}}\Big\|\sum_{n \in I_+}c_nh_n\Big\|.\end{eqnarray*}Similarly, as $\{h_n\}_{n \in I_-}$ is complete in $\mathcal{K}^-_1$  and by taking $\{c_n\}_{n \in I_+}\in \ell_2(I_-)$, we have
		\begin{eqnarray*}
			\Big\|U^-\Big(\sum_{n \in I_-}c_nh_n\Big)\Big\| &\leq&\sqrt{\frac{B^{'}}{A^{'}}}\Big\|\sum_{n \in I_-}c_nh_n\Big\|.\end{eqnarray*} Hence $U^+$ and $U^-$ are bounded operators and have extensions to bounded operators on $\mathcal{K^+}$ and $\mathcal{K^+}$ respectively. 	The following theorem gives a  characterization for Riesz bases in Krein spaces.
	\begin{theorem}\label{rt} Let $\{f_n\}_{n \in I}$ be a sequence in $(\mathcal{K},[.,.]).$ Then the following statements are equivalent:
		\begin{enumerate}
			\item [(a)] $\{f_n\}_{n \in I}$ is a Riesz basis for $(\mathcal{K},[.,.]).$
			\item [(b)] $\{f_n\}_{n \in I}$ is complete in $(\mathcal{K},[.,.])$ and there exist positive constants $A,B,A^{'}$ and $B^{'}$ such that for $\{c_n\}_{n \in I_+} \in \ell_2(I_+),$ then 
			\begin{equation}\label{a1}
			A\sum_{n \in I_+}|c_n|^2 \leq \Big\|\sum_{n \in I_+}c_nf_n\Big\|^2 \leq B\sum_{n \in I_+}|c_n|^2
			\end{equation}   
			and for $\{c_n\}_{n \in I_-} \in \ell_2(I_-)$ we have
			\begin{equation}\label{a2}
			A^{'}\sum_{n \in I_-}|c_n|^2 \leq \Big\|\sum_{n \in I_-}c_nf_n\Big\|^2 \leq B^{'}\sum_{n \in I_-}|c_n|^2.
			\end{equation} 
		\end{enumerate}
	\end{theorem}
	\begin{proof}(a) $\Rightarrow$ (b) :
		First suppose that $\{f_n\}_{n \in I}$ is a Riesz basis for the Krein space $(\mathcal{K},[.,.]),$ write $f_n= \{\{U^+e_n\}_{n \in I_+} \cup \{U^-e_n\}_{n \in I_-}\}$ by definition. Also $\{f_n\}_{n \in I_+}$ is complete by the consequence of the Theorem \ref{T1} and for any $\{c_n\}_{n \in I_+} \in \ell_2(I_+)$ we have 
		\begin{eqnarray*}
			\Big\|\sum_{n \in I_+}c_nf_n\Big\|^2&=& \Big\|U^+\Big(\sum_{n \in I_+}c_ne_n\Big)\Big\|^2 \\&\leq&\|U^+\|^2\Big\|\sum_{n \in I_+}c_ne_n\Big\|^2\\&=& \|U^+\|^2\sum_{n \in I_+}|c_n|^2
		\end{eqnarray*} and \begin{eqnarray*}
			\Big\|\sum_{n \in I_+}c_ne_n\Big\|^2&=&\Big\|{U^+}^{-1}U^+\Big(\sum_{n \in I_+}c_ne_n\Big)\Big\|^2
			\\&\leq&\|{U^+}^{-1}\|^2\Big\|\sum_{n \in I_+}c_nf_n\Big\|^2.
		\end{eqnarray*}
	Hence
		\begin{eqnarray*}
			\frac{1}{\|{U^+}^{-1}\|^2}\sum_{n \in I_+}|c_n|^2 \leq \Big\|\sum_{n \in I_+}c_nf_n\Big\|^2 \leq \|U^+\|^2\sum_{n \in I_+}|c_n|^2.
		\end{eqnarray*} Similarly, by Theorem \ref{T1}, $\{f_n\}_{n \in I_-}$ is complete and for $\{c_n\}_{n \in I_+} \in \ell_2(I_-)$ we get 
		\begin{eqnarray*}
			\frac{1}{ \|{U^-}^{-1}\|^2}\sum_{n \in I_-}|c_n|^2 \leq \Big\|\sum_{n \in I_-}c_nf_n\Big\|^2 \leq \|U^-\|^2\sum_{n \in I_-}|c_n|^2.\end{eqnarray*}
		(b)$\Rightarrow$ (a) : The inequalities  (\ref{a1}) and (\ref{a2}) imply that $\{f_n\}_{n \in I_+}$ and $\{f_n\}_{n \in I_-}$ are Bessel sequences with Bessel bounds $B$ and $B^{'}$ respectively. Let $\{e_n\}_{n \in I_+}$ and $\{e_n\}_{n \in I_-}$ be orthonormal bases for $\mathcal{K^+}$ and $\mathcal{K^-}$ respectively and we can extend $U^+_1e_n= f_n$ to a bounded operators on   $\mathcal{K}^+_1$ and $U^-_1e_n= -f_n$ to a bounded operator on $\mathcal{K}^-_1.$ In a similar fashion, extend ${U_2}^+f_n=e_n$ to a bounded operator on $\mathcal{K}^+_2$ and ${U_2}^-f_n=-e_n$ to a bounded operator on $\mathcal{K}^-_2.$ Then ${U_2}^+U^+_1=U^+_1{U_2}^+=I$ and ${U_2}^-U^-_1=U^-_2{U_2}^-=I,$ so $U^+_1$ and $U^-_1$ are invertible operators. Hence $\{f_n\}_{n \in I}$ is a Riesz basis for $(\mathcal{K},[.,.]).$ 
	\end{proof}
	\section{Gram matrix} 
	In this section, we introduce the concept of Gram matrix in $(\mathcal{K},[.,.])$ and prove some of its properties. Let $\{f_n\}_{n \in I}\subset \mathcal{K}$ be a Bessel sequence, we can compose the synthesis operator $T$ and its adjoint $T^*$; where $T=T^++T^-$ and $T^*={T^+}^*+{T^-}^*.$ We obtain the bounded operators ${T^+}^*T^+: \ell_2(I_+)\longrightarrow \ell_2(I_+)$ defined by
	\begin{eqnarray*}
		{T^+}^*T^+\{c_n\}_{n \in I_+}=\Big\{[\sum_{n \in I_+}c_nf_n,f_n]\Big\}
	\end{eqnarray*} and ${T^-}^*T^-: \ell_2(I_-)\longrightarrow \ell_2(I_-)$ defined by
	\begin{eqnarray*}
		{T^-}^*T^-\{c_n\}_{n \in I_-}=\Big\{[\sum_{n \in I_-}c_nf_n,f_n]\Big\}
	\end{eqnarray*}
	Let $\{e_n\}_{n \in I_+}$ and $\{e_n\}_{n \in I_-}$ be the canonical orthonormal bases for $\ell_2(I_+)$ and $\ell_2(I-)$ respectively.  Then the $ij$-th entries in the matrix representations for ${T^+}^*T^+$ and ${T^-}^*T^-$ are
	\begin{eqnarray*}
		[ {T^+}^*T^+ e_i,e_j ]=[T^+ e_i,T^+e_j ]\Big\{([f_i,f_j]_{i,j \in I_+}\Big\}
	\end{eqnarray*} and 
	\begin{eqnarray*}
		[ {T^-}^*T^- e_i,e_j ]=[T^- e_i,T^-e_j ]=\Big\{[f_i,f_j]_{i,j \in I_-}\Big\}.
	\end{eqnarray*}
	The matrix $\Big\{[f_i,f_j]_{i,j \in I_+}\Big\}$ is called the positive Gram matrix and $\Big\{[f_i,f_j]_{i,j \in I_-}\Big\}$ is called the negative Gram matrix associated with $\{f_n\}_{n \in I_+}$ and $\{f_n\}_{n \in I_-}$ respectively. Moreover, if $\{f_n\}_{n \in I}\subset \mathcal{K}$ is a Bessel sequence, then the Gram matrix defines bounded operators on $\ell_2(I_+)$ and $\ell_2(I_-).$\\
	The following lemma states that in order to prove that the Gram matrices are bounded we can consider the Bessel condition.
	\begin{lemma}\label{lg}Let $\{f_n\}_{n \in I}$ be a sequence in $(\mathcal{K},[.,.])$. Then the following conditions are equivalent:
		\begin{enumerate}
			\item [(a)] $\{f_n\}_{n \in I}$ is a Bessel sequence with Bessel bound $B$.
			\item [(b)] The positive Gram matrix and the negative Gram matrix  associated with $\{f_n\}_{n \in I_+}$ and $\{f_n\}_{n \in I_-}$ define bounded operators on $\ell_2(I_+)$ and $\ell_2(I_-),$ with norm at most $B$ and $B^{'}$ respectively.
		\end{enumerate}
	\end{lemma}
	\begin{proof}$(a)\Rightarrow (b)$ : Suppose $\{f_n\}_{n \in I}$ is a Bessel sequence in $(\mathcal{K},[.,.])$. Then  the positive Gram matrix and the negative Gram matrix define bounded operators on $\ell_2(I_+)$ and $\ell_2(I-)$ and the proof for the norm estimation follows from the Theorem 3.1 in \cite{s}.
		 
		$(b)\Rightarrow (a)$ :  Next, suppose that $(b)$ holds. Let $\{c_n\}_{n \in I_-} \in \ell_2(I_-).$ Then
		\begin{eqnarray*}
			\sum_{j \in I_-}\Big|\sum_{n \in I_-}c_n[f_n,f_j]\Big|^2\leq B^2\sum_{n \in I_-}|c_n|^2.
		\end{eqnarray*}
		Let $\ell,m \in I_-$ such that $\ell>m.$ Then
		\begin{eqnarray*}
			\Big\|\sum_{n=1}^\ell c_nf_n-\sum_{n=1}^mc_nf_n\Big\|^4 &=& \Big\|\sum_{n =m+1}^\ell c_nf_n\Big\|^4
			\\&=&\Big|[\sum_{n=m+1}^\ell c_nf_n,\sum_{k=m+1}^\ell c_kf_k]\Big|^2
			\\&=&\Big|\sum_{n=m+1}^\ell \overline{c_n} \sum_{k=m+1}^\ell c_k[f_n,f_k]\Big|^2
			\\&\leq&\Big(\sum_{n=m+1}^\ell|c_n|^2\Big)\Big(\sum_{n=m+1}^\ell\Big|\sum_{k=m+1}^\ell c_k[f_n,f_k]\Big|^2\Big).
		\end{eqnarray*}
		Therefore 
		\begin{eqnarray*}
			\sum_{n=m+1}^\ell\Big|\sum_{k=m+1}^\ell c_k[f_n,f_k]\Big|^2 \leq B^2\sum_{n=m+1}^\ell|c_n|^2.
		\end{eqnarray*}
Hence we conclude that \begin{eqnarray*}
			\Big\|\sum_{n=1}^\ell c_nf_n-\sum_{n=1}^mc_nf_n\Big\|^4\leq B^2\Big(\sum_{n=m+1}^\ell |c_n|^2\Big)^2.
		\end{eqnarray*}
		Thus $\sum_{n \in I_-}c_nf_n$ is convergent and by repeating the above argument, we get
		\begin{eqnarray*}
			\Big\|\sum_{n \in I_-}c_nf_n\Big\| \leq \sqrt{B}\Big(\sum_{n \in I_-} |c_n|^2\Big)^\frac{1}{2}.
		\end{eqnarray*} Similarly, one can easily prove that 
		\begin{eqnarray*}
			\Big\|\sum_{n \in I_+}c_nf_n \Big\| \leq \sqrt{B^{'}}\Big(\sum_{n \in I_+}|c_n|^2\Big)^\frac{1}{2}.
		\end{eqnarray*}
		Hence by Theorem 3.1 in \cite{s}, $\{f_n\}_{n \in I}\subset \mathcal{K}$ is a Bessel sequence with Bessel bound B.
	\end{proof} 
The following lemma gives us a sufficient condition for the Gram matrices defining bounded operators on $\ell_2(I_+)$ and $\ell_2(I_-).$
	\begin{lemma} Let $M_1=\{{M_1}_{j,n}\}_{j,n \in I_+}$ and let $M_2=\{{M_2}_{j,n}\}_{j,n \in I_-}$ be two matrices for which ${M_1}_{j,n} =\overline{{M_1}_{n,j}}$ for all $j,n \in I_+$ and ${M_2}_{j,n} =\overline{{M_2}_{n,j}}$ for all $j,n \in I_-$ and  there exist constants $B$ and $B^{'}$ such that $\sum_{j,n \in I_+}|{M_1}_{j,n}| \leq B$ and $\sum_{j,n \in I_+}|{M_2}_{j,n}| \leq B^{'}.$ Then $M_1$ and $M_2$ define bounded linear operators on $\ell_2(I_+)$ and $\ell_2(I_-)$ with norm at most $B$ and $B^{'}$ respectively. 
	\end{lemma}
	\begin{proof}
		Let $\{c_n\}_{n \in I_+} \in \ell_2(I_+).$ Then ${M_1}\{c_n\}_{n \in I_+}$ is well defined as the sequence indexed by $I_+$ and whose $j$-th coordinate is given by $\sum_{n \in I_+}{M_1}_{j,n}c_n.$ To prove this sequence is in $\ell_2(I_+),$  it is enough to show that the map \begin{eqnarray}\label{lm}\phi:\{d_n\}_{n \in I_+} \longrightarrow [\{d_n\}_{n \in I_+},{M_1}\{c_n\}_{n \in I_+}]_{\ell_2(I_+)} \end{eqnarray} is a continuous linear functional on $\ell_2(I_+).$ Now, for $\{d_n\}_{n \in I_+} \in \ell_2(I_+),$ we have 
		\begin{eqnarray*}
			\sum_{j \in I_+}\Big|\sum_{n \in I_+}\overline{{M_1}_{n,j}c_n}d_j \Big| &\leq& \sum_{j \in I_+}\sum_{n \in I_+}|{M_1}_{j,n}c_nd_j|
			\\&=& \sum_{j \in I_+}\sum_{n \in I_+}\Big(|{M_1}_{j,n}|^\frac{1}{2}|c_n|)(|{M_1}_{j,n}|^\frac{1}{2}|d_j|\Big)
			\\&\leq& \Big(\sum_{j \in I_+}\sum_{n \in I_+}|{M_1}_{j,n}||c_n|^2\Big)^\frac{1}{2}\Big(\sum_{j \in I_+}\sum_{n \in I_+}|{M_1}_{j,n}||d_j|^2\Big)^\frac{1}{2}
			\\&\leq& B\Big(\sum_{n \in I_+}|c_n|^2\Big)^\frac{1}{2}\Big(\sum_{j \in I_+}|d_j|^2\Big)^\frac{1}{2}.
		\end{eqnarray*}From above, we conclude that $(\ref{lm})$ defines a continuous linear functional on $\ell_2(I_+).$ Also,
		\begin{eqnarray*}
			\|{M_1}\{c_n\}_{n \in I_+}\| &=&\sup\limits_{\|\{d_n\}\|=1}\Big|[\{d_n\}_{n \in I_+}, {M_1}\{c_n\}_{n \in I_+}]_{\ell_2(I_+)}\Big|
			\\&\leq& B\Big(\sum_{n \in I_+}|c_n|^2\Big)^\frac{1}{2}.
		\end{eqnarray*}
		This shows that ${M_1}$ is a bounded linear operator on $\ell_2(I_+)$ with norm at most $B.$ Similarly, it is easy to prove that ${M_2}$ defines a bounded linear operator on $\ell_2(I_+)$ with norm at most $B^{'}.$ 
	\end{proof}
	\begin{proposition}Let $\{f_n\}_{n \in I}$ be a sequence in  $(\mathcal{K},[.,.])$.  If there exists a positive constant $B$ such that 
		\begin{eqnarray}\label{pm}
		\sum_{j, n \in I}|[f_j,f_n]| \leq B,
		\end{eqnarray} then $\{f_n\}_{n \in I}$ is a Bessel sequence with Bessel bound B.  
	\end{proposition} The condition $(\ref{pm})$ is much better than that of Bessel condition in \cite{s} (Equation 6). Because here we have to check the indefinite inner product between the elements of $\{f_n\}_{n \in I},$ whereas in the case of Bessel condition we have to check the indefinite inner product for all $f$ in the Krein space $\mathcal{K}$.

	The following theorem gives an equivalent condition of the Riesz basis in terms of Gram matrices.
	\begin{theorem}
		Let $\{f_n\}_{n \in I}$ be a sequence in  $(\mathcal{K},[.,.])$. Then the following conditions are equivalent:
		\begin{enumerate}
			\item [$(a)$] $\{f_n\}_{n \in I}$ is a Riesz basis for  $(\mathcal{K},[.,.]).$
			\item [$(b)$] $\{f_n\}_{n \in I_+}$ and $\{f_n\}_{n \in I_-}$ are complete and their Gram matrices $\{[f_n,f_j]_{n,j \in I_+}\}$ and $\{[f_n,f_j]_{n,j \in I_-}\}$ define bounded, invertible operators on $\ell_2(I_+)$ and $\ell_2(I_-)$ respectively.
	\end{enumerate}\end{theorem}
	\begin{proof} (a) $\Rightarrow$ (b) :
		Since $\{f_n\}_{n \in I}$ is a Riesz basis, we can write $\{f_n\}_{n \in I}= \{\{U^+e_n\}_{n \in I_+}\cup \{U^-e_n\}_{n \in I_-}\}$ by definition. For any $n,j \in I_+,$ we have 
		\begin{eqnarray*}
			[f_n,f_j] = [U^+e_n,U^+e_n]= [{U^+}^*U^+e_n,e_j]
		\end{eqnarray*} and for $n,j \in I_-$
		\begin{eqnarray*}
			[f_n,f_j]=[{U^-}^*U^-e_n,e_j].
		\end{eqnarray*}
		It is the positive Gram matrix  representing the bounded invertible operator ${U^+}^*U^+$ in the basis $\{e_n\}_{n \in I_+}$ and the negative Gram matrix is the matrix representing the bounded invertible operator ${U^-}^*U^-$ in the basis $\{e_n\}_{n \in I_-}.$
		
		$(b)\Rightarrow (a)$ : Suppose that (b) holds.  Then the upper Riesz conditions in $(\ref{a1})$ and $(\ref{a2})$ are satisfied by using Lemma (\ref{lg}) and Theorem $(\ref{rt}).$ Next, we  will prove that the lower conditions in $(\ref{a1})$ and $(\ref{a2})$ are also satisfied. For this, let $G_1$ denote the operator on $\ell_2(I_+)$ given by the Gram matrix $\{[f_n,f_j]_{n,j \in I_+}\}$ and $G_2$ denote the operator on $\ell_2(I_-)$ given by the Gram matrix $\{[f_n,f_j]_{n,j \in I_-}\}.$ 
		
		First, for a given sequence $\{c_n\}_{n \in I_+} \in \ell_2(I_+)$, the $j^{\text{th}}$ element in the image sequence $G_1\{c_n\}$ is $\sum_{n \in I_+}[f_n,f_j]e_n.$ So,
		\begin{equation*}
			[G_1 \{c_n\}_{n \in I_+}, \{c_k\}_{k \in I_+}] = \sum_{j \in I_+}\sum_{k \in I_+}[f_n,f_j]c_n\overline{c_j}
			= \|\sum_{n \in I_+}c_nf_n\|.
		\end{equation*} 
		Thus $G_1$ is positive and with a similar argument it is easy to show that $G_1$ is self adjoint. Let $W$ be a positive square root of $G_1$. Then by the above calculation, we get 
		\begin{eqnarray*}
			\Big\|\sum_{n \in I_+}c_nf_n\Big\| &=&\|W\{c_n\}_{n \in I_+}\|^2 \ \ \geq \ \  \frac{1}{\|W^{-1}\|^2}\sum_{n \in I_+}|c_n|^2.
		\end{eqnarray*}
	\end{proof}

\begin{center}
	\textbf{Acknowledgements}
\end{center}

\noindent The first author is highly thankful to the Central University of Haryana for providing basic facilities to carry out this research.  The present work of the second author is partially supported by  ANRF (SERB), DST, Government of India (TAR/2022/000219) under the scheme “TARE”.

\end{document}